\documentclass[12pt]{article}
\usepackage{epsfig}

\newtheorem{lem}{Lemma}
\newtheorem{thm}{Theorem}

\newtheorem{rem}{Remark}

\begin{document}
\centerline{\large\bf On bifurcations of multidimensional diffeomorphisms having}
\centerline{\large\bf a homoclinic tangency to a saddle-node.}
\medskip
\centerline{\bf Gonchenko S.V.$^1$, Gordeeva O.V.$^1$, Lukjanov V.I.$^1$, Ovsyannikov I.I.$^{12}$}
\medskip
\centerline{\it $^1$ Nizhny Novgorod State University}
\centerline{\it $^2$ Universit\"at Bremen, Germany}

\centerline{\it E-mail: gonchenko@pochta.ru, olga.gordeeva@inbox.ru,}
\centerline{\it viluk@yandex.ru, Ivan.I.Ovsyannikov@gmail.com}
\bigskip

\abstract{We study main bifurcations of multidimensional diffeomorphisms having a non-transversal homoclinic
orbit to a saddle-node fixed point. On a parameter plane we build a bifurcation diagram for single-round periodic orbits
lying entirely in a small neighbourhood of the homoclinic orbit. Also a relation of our results to
well-known codimension one bifurcations of a saddle fixed point with a quadratic homoclinic tangency and
a saddle-node fixed point with a transversal homoclinic orbit is discussed.
}



\section{Introduction}
Bifurcations of systems with homoclinic structures have a special meaning for the mathematical theory of dynamical chaos. It is well-known from the Shilnikov work \cite{Sh67} that  the set of orbits entirely lying in a neighbourhood of a transverse Poincar\'e homoclinic orbit has a complex structure: it is the nontrivial uniformly hyperbolic set containing a countable number of periodic and hyperbolic orbits, continuum Poisson stable orbits etc. In the case of systems with homoclinic tangencies, the situation becomes much more complicated and even unpredictable in a sense. The point is that bifurcations of such systems can lead to the appearance of periodic and homoclinic orbits of any orders of degeneracy, \cite{GTS91, GST93, GTS99, GST07}. Therefore, the complete study of bifurcations of such systems is principally impossible and, hence, the problems connecting with the study of principal bifurcations and characteristic features of the dynamics should come to the foreground here. This fully relates to the problem under consideration
in the present paper ---  the study of bifurcations of
diffeomorphisms with homoclinic tangencies to saddle-node fixed points.

\begin{figure}[ht]
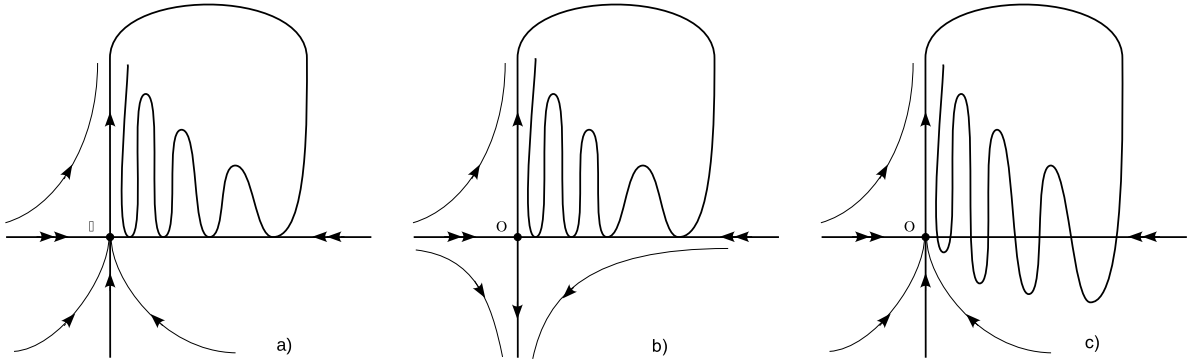

\begin{center}
\begin{tabular}{ccc}
\epsfig{file=initial.eps,width=5cm}  &
\epsfig{file=saddle.eps,width=5cm} &
\epsfig{file=intersect.eps,width=5cm}
\end{tabular}
\caption{{\footnotesize a) saddle-node with a homoclinic tangency;
b) saddle with a homoclinic tangency;
c) saddle-node with a transversal homoclinic orbit.}}
\label{fig:options}
\end{center}
\end{figure}

This problem is naturally connected with two other known problems concerning the study of
global bifurcations. The first of them  is
a bifurcation of diffeomorphisms with quadratic homoclinic tangencies to a saddle fixed point. Note that it
was first studied
by N.K. Gavrilov and L.P. Shilnikov in their famous paper \cite{GaS73}, where, in fact, the foundation of the mathematical theory of homoclinic chaos was laid.
The second problem, the study of global bifurcations of diffeomorphisms with a transverse homoclinic orbit to a saddle-node, was considered first by V.I. Lukjanov and L.P. Shilnikov in  \cite{LukSh78} in which certain
conditions for the
onset of  chaos immediately after the disappearance of the saddle-node were found. In fact, in  \cite{LukSh78} the mathematical theory was constructed that has predicted and explained the phenomenon ``intermittency'' known as one of main mechanism of chaos onset in models from applications.

The main results of the present paper, Theorems 1 and 2, can be considered as a generalization of  the well-known theorems on ``cascade of periodic sinks'',  \cite{GaS73, LukSh78, G83},  onto the case of homoclinic tangency to a saddle-node. One can say that our case (Fig.~\ref{fig:options}(a))
is a ``meeting point'' of the  Gavrilov-Shilnikov (Fig.~\ref{fig:options}(b)) and Lukjanov-Shilnikov  
(Fig.~\ref{fig:options}(c)) cases. This accounts for our interest to the problem under consideration.

We note that the study of
global bifurcations accompanying
the disappearance of a saddle-node with homoclinic orbits 
has also an important meaning
for applications. 
In particular, such bifurcations underly certain scenarios of the emergence of strange attractors of type ``torus-chaos''. 
Physically, these scenarios describe bifurcation phenomena which occur  at the transition
from a synchronization regime to a chaotic one. Mathematical basics of the corresponding theory were laid in the work by V.S. Afraimovich and L.P. Shilnikov \cite{AfSh83}, see also \cite{NPT,ST00}.

\begin{figure}[ht]
\begin{center}
\epsfig{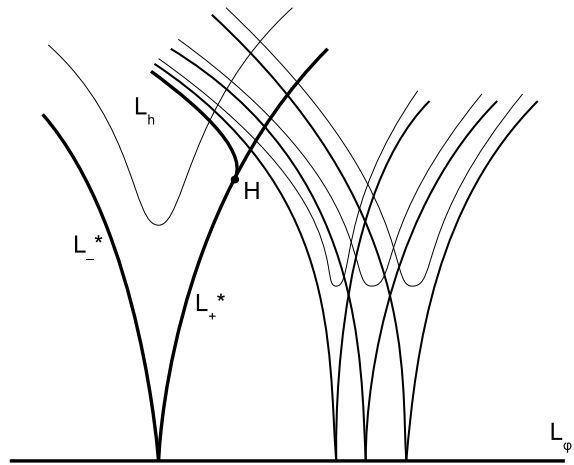}
\vspace{-1cm}
\caption{{\footnotesize The Arnold tongues.}}
\label{fig:tongue}
\end{center}
\end{figure}

\begin{figure}[ht]
\begin{center}
\epsfig{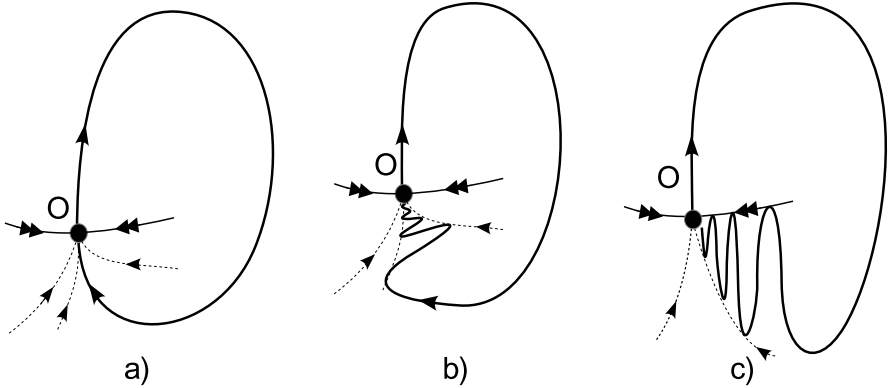}
\vspace{-1cm}
\caption{{\footnotesize a) smooth; b) non-smooth closed invariant curve; c) the homoclinic tangency appearance.}}
\label{fig:curve}
\end{center}
\end{figure}

In these scenarios, the codimension two bifurcation of homoclinic tangency to a saddle-node periodic point plays a significant role. Thus, in Fig.~\ref{fig:tongue}
a fragment of the bifurcation diagram near the
line $L_\varphi$ corresponding to the
existence of a fixed point with multipliers $e^{\pm i\varphi}$ is shown. When parameters cross $L_\varphi$ (upwards) the Neimark-Sacker bifurcation occurs and, as a result, the fixed point becomes unstable (of focal type) and a stable closed invariant curve is born in its neighbourhood. It is well-known that a pair of bifurcation lines $L_{p/q}^1$ and $L_{p/q}^2$ starts from every resonant point $\varphi = 2\pi p/q$ of $L_\varphi$ corresponding to the existence of a saddle-node point of period $q$ on the invariant curve. 
 The domain between the lines $L_{p/q}^1$ and $L_{p/q}^2$ is called 
a synchronization zone or an ``Arnold tongue''. The saddle-node falls into two points (saddle and node) of period $q$ for values of parameters inside the tongue and disappears outside it.  When changing values of parameters along the line $L_{p/q}^{1,2}$ the invariant curve is  smooth at the beginning, Fig.~\ref{fig:curve}(a), then it loses the smoothness and becomes ``corrugated'',  Fig.~\ref{fig:curve}(b). Moreover, a point $H$ exists on $L_{p/q}^{1,2}$ such that  the unstable and strong stable invariant manifolds of the saddle-node touch,  Fig.~\ref{fig:curve}(a), evidently, at this moment the invariant curve does not longer exist. Accordingly, a bifurcation line $L_h$ exists inside the tongue  beginning at the point $H$ and corresponding to the existence of a homoclinic tangency to the saddle point. 

Note that the study of bifurcations of diffeomorphisms in a
neighbourhood of the point $H$ has never been done before and, thus, it remains an important stage for
understanding bifurcation phenomena that occur during the transition from  synchronization regime to chaos. We need to recall several 
results related to this topic.
In \cite{Gor07a}, \cite{Gor07b} a structure of orbits entirely lying in a small
neighbourhood of homoclinic orbits to nonhyperbolic points (of saddle-node or saddle type) was studied. Bifurcations of 
two-dimensional flows having a homoclinic loop $\Gamma$ of a saddle-node equilibrium $O$ such that $\Gamma \subset W^{ss}$ (i.e. $\Gamma$ enters $O$ along the nonleading  direction) were considered in \cite{Luk82}. 
Note that homoclinic tangencies to  saddle-nodes  can naturally appear under small periodic perturbations of such  autonomous flows.

In the present paper we study principal bifurcations in two parameter families of
multidimensional diffeomorphisms having a quadratic homoclinic tangency to a saddle-node.

\section {Statement of the problem and main results}
Consider a $C^r$-smooth, $r \ge 4$, $(n + 1)$-dimensional diffeomorphism $f_0$, $n \ge 1$, having a fixed
point $O$ of a saddle-node type and a nontransversal homoclinic to it orbit $\Gamma_0$,
see fig.~\ref{fig:options}~a). We assume that $f_0$ has no other degeneracies
i.e. it satisfies conditions {\bf A.} and {\bf B.} formulated below.

{\bf A.} The fixed point $O$ is a non-degenerate saddle-node, i.e. it has multipliers
$\lambda_1, \lambda_2, \ldots, \lambda_{n+1}$, where $|\lambda_i| < 1$, $i = \overline{1, n}$, and $\lambda_{n + 1} = 1$.
The first Lyapunov value $l_1$ at $O$ is non-zero. Without loss of generality we assume that $l_1 > 0$.

Let $U_0$ be some small neighbourhood of $O$. It is well known that there exists a nonleading (strong stable)
$C^r$-smooth $n$-dimensional invariant manifold $W^{ss}(O)$ tangent to the eigendirections corresponding to
multipliers $\lambda_1, \lambda_2, \ldots, \lambda_{n}$. Thus, $U_0$ is divided by $W^{ss}$ into two parts
which we will call the node $U^+$ and the saddle $U^-$ zones. The positive semi-trajectory of any point of $U^+$
tends to $O$ being tangent to the leading eigendirection corresponding to the multiplier $\lambda_{n + 1} = 1$. In the saddle
zone $U^-$ there exists a one-dimensional $C^r$-smooth unstable manifold $W^u$ which consists of those points
negative semi-trajectories of which tend to $O$ at backward iterations of $f_0$. All other points of $U^-$ leave $U_0$
at the both forward and backward iterations. The unstable manifold $W^u$ is tangent to the leading eigendirection
(corresponding to $\lambda_{n + 1}$) at $O$.

{\bf B.} Invariant manifolds $W^u(O)$ and $W^{ss}(O)$ have a quadratic tangency at the points of the homoclinic
orbit $\Gamma_0$.

Take some homoclinic point $M^+ \in W^{ss} \cap U_0$ of $\Gamma_0$ and its small neighbourhood $\Pi^+ \subset U_0$.
Denote as $l_u$ a piece $W^u \cap \Pi^+$ of the unstable manifold $W^u$ containing $M^+$. We will distinguish two main cases
of a homoclinic tangency to the saddle-node, {\em tangency from a node zone} when $l_u \subset U^-$ (fig.~\ref{fig:options}a)
and {\em tangency from a saddle zone} when $l_u \subset U^+$ (fig.~\ref{fig:curve}c).

We recall the following facts from the theory of invariant manifolds \cite{ST00, Luk82, HPS, book}.
First of all, $f_0$ in $U_0$ possesses a
center manifold $W^C(O)$ which is not unique. Each such manifold is a one-dimensional $C^r$-smooth curve
which coincides with $W^u_{loc}(O)$ in $U^-$, intersects $U^+$ and is tangent in it to the leading eigendirection.
Moreover, $f_0$ in $U_0$ possesses a unique strong stable invariant foliation $F^{ss}$. Its leaves are $n$-dimensional planes
which are transversal to each $W^C(O)$, in particular, in $U^-$ they all are transversal to $W^u$. Moreover, the strong stable
manifold $W^{ss}$ is one of the leaves.

Diffeomorphisms satisfying  conditions {\bf A} and {\bf B} comprise a codimension-2
bifurcation surface ${\cal B}_2$ in the space of $(n+1)$-dimensional $C^r$-diffeomorphisms.
To study bifurcations of $f_0$ we consider a two parameter
family $f_\mu$, $\mu = (\mu_1, \mu_2)$, which is transverse to ${\cal B}_2$ at $\mu = 0$.
The initial diffeomorphism $f_0$ belongs to this family at $\mu = 0$.

We choose governing parameters $\mu_1$ and $\mu_2$ in the following way.
Parameter $\mu_1$ is the splitting parameter of manifolds $W^u$ and $W^{ss}$ with respect to some homoclinic point
(for example, $M^+$) in the case when the saddle-node exists. Thus, $W^u$ and $W^{ss}$ will intersect transversely
in two points close to $M^+$ when $\mu_1 < 0$ in the case of tangency from the saddle zone ($\mu_1 > 0$ in the case
of tangency from the node zone); correspondingly, $W^u$ and $W^{ss}$ will have no intersection near $M^+$ if $\mu_1 > 0$
($\mu_1 < 0$). The second parameter $\mu_2$ controls bifurcations of the fixed point $O$ such that
$f_\mu$ possesses a saddle-node fixed point at $\mu_2 = 0$, does not have fixed points in $U_0$ (the saddle-node disappears)
when $\mu_2 > 0$ and the saddle-node is split into two fixed points (saddle $O_1$ and node $O_2$) when $\mu_2 < 0$.

While changing parameters $(\mu_1, \mu_2)$, bifurcations will take place in $f_\mu$. In particular, they will be
related to the emergence of new homoclinic structures. 

%
%
\begin{figure}[ht]
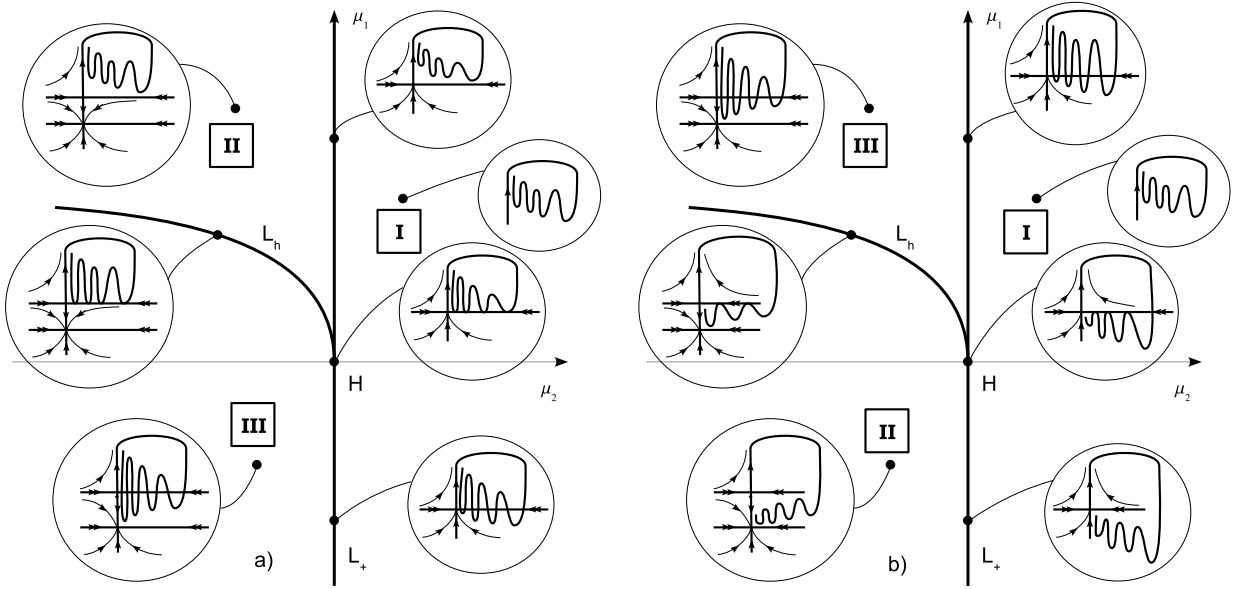

\begin{center}
\begin{tabular}{cc}
\epsfig{file=BifDiagrLuk.eps,width=8cm, height=8cm}  &
\epsfig{file=BifDiagrLukbelow.eps,width=8cm, height=8cm}
\end{tabular}
\caption{Splitting the $(\mu_1, \mu_2)$ plane onto domains by the codimension-one bifurcation lines $L_+$ and $L_h$.
a) Tangency from the saddle zone; b) tangency from the node zone.}
\label{fig:bifdiagluk}
\end{center}
\end{figure}
%
%
%

%
\begin{thm}
\label{thm:bifdiag}
In any sufficiently small neighbourhood of the origin $\mu = 0$ on the parameter plane $\mu = (\mu_1, \mu_2)$ 
there exist two bifurcation
curves $L_+: \; \mu_2 = 0$ and $L_h: \; \mu_1 = \sqrt{-\mu_2} + O(|\mu_2|^{3/2})$, see Fig.~\ref{fig:bifdiagluk}.
If $\mu \in L_+$, each diffeomorphism $f_\mu$ possesses a fixed point of the saddle-node type close to $O$; if
$\mu \in L_h$, then  $f_\mu$ has a nontransversal homoclinic orbit $\Gamma_\mu$ close to $\Gamma_0$ in points of which the
invariant manifolds $W^u(O_1)$ and $W^s(O_1)$ of the saddle fixed point $O_1$ have a quadratic tangency.

\begin{itemize}
\item In domain {\bf I} $(\mu_2 > 0)$ the diffeomorphism $f_\mu$ has no fixed points.

\item In domains {\bf II} and {\bf III} the diffeomorphism $f_\mu$ has two fixed points, saddle $O_1$ and node $O_2$. If $\mu \in \mbox{\bf II}$,
there are no homoclinic orbits to saddle $O_1$ close to $\Gamma_\mu$, and if $\mu \in \mbox{\bf III}$,
there are two transversal homoclinic orbits to $O_1$ close to $\Gamma_\mu$.

\end{itemize}

\end{thm}

Figures~\ref{fig:bifdiagluk} illustrate the statement of the theorem in both cases of
tangency from the saddle and node zones respectively.

Consider a sufficiently small neighbourhood $U$ of $O \cup \Gamma_0$ which is called an extended
neighbourhood of the saddle-node $O$. It is a union
of a small neighbourhood $U_0$ of the fixed point $O$ and a finite set of
small neighbourhoods $V_i$ of those points of $\Gamma_0$ that lie outside $U_0$.
We call a periodic orbit of $f_\mu$ lying entirely in $U$ {\em single-round} if it intersects
each of $V_i$ in exactly one point.

We will study bifurcations of single-round periodic orbits in the two parameter family $f_\mu$ for small $\mu$.

For $\mu = 0$ we choose in $U_0$ two homoclinic points $M^+ \in W^{ss}_{loc}(O)$ and
$M^- \in W^{u}_{loc}(O)$. Let $\Pi^+, \Pi^- \subset U_0$ be their respective sufficiently small neighbourhoods
of diameter $\varepsilon$. Evidently,
there exists an integer $q$ such that $M^- = f_0^q(M^+)$. Then,
for all small $\mu$,  two maps by the orbits of the diffeomorphism $f_\mu$ will be defined in $U_0$: these are
the {\em local map} $T_0 = \left. f_\mu \right|_{U_0}$
and the {\em global map} $T_1 = f_\mu^q:\; \Pi^- \to \Pi^+$. By construction, any single-round periodic orbit, lying entirely in $U$,
has exactly one intersection point with the neighbourhoods $\Pi^+$ and $\Pi^-$. A point of such an orbit in $\Pi^+$ can be
regarded as a fixed point of the corresponding {\em first return map} $T_k = T_1 \cdot T_0^k: \; \Pi^+ \to \Pi^- \to \Pi^+$,
where $T_0^k$ maps points from $\Pi^+$ to $\Pi^-$. We will study bifurcations of fixed points of the first return maps
$T_k$ for all sufficiently large $k$. As a result we will construct the bifurcation diagram for single-round periodic
orbits of diffeomorphisms from the family $f_\mu$. Its structure is described by the following theorem:


\begin{figure}[ht]
\begin{center}
\epsfig{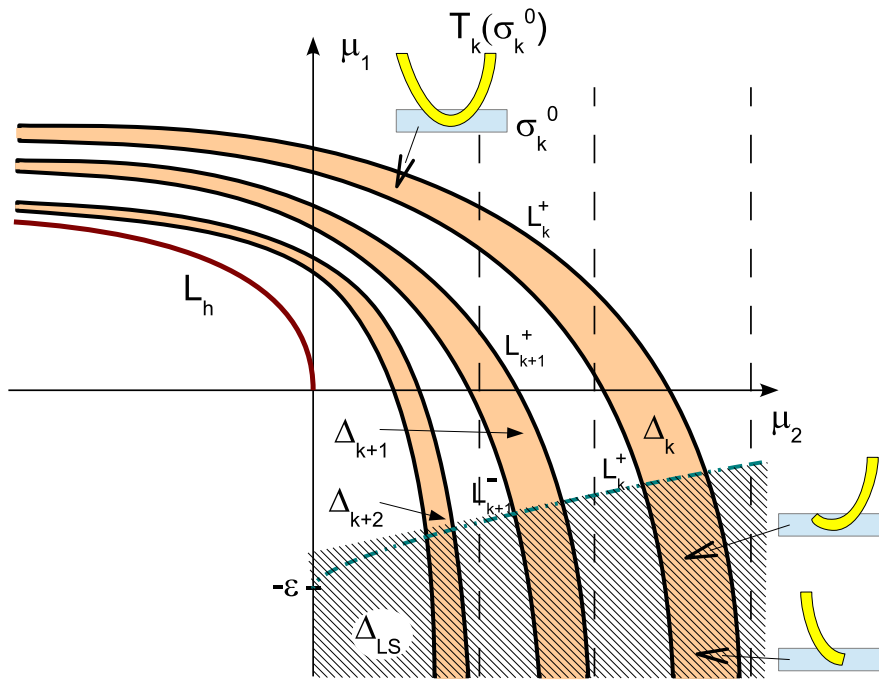}
\vspace{-1cm}
\caption{{\footnotesize Bifurcation diagram for the family $f_\mu$ in the case of tangency from the saddle zone.}}
\label{fig:bifdiag}
\end{center}
\end{figure}

\begin{thm}
\label{thm:main}
1$)$ On the parameter plane $(\mu_1, \mu_2)$ in any sufficiently small neighbourhood of the origin there exist
a countable set of disjoint domains $\triangle_k$ such that for $\mu \in \triangle_k$ the diffeomorphism $f_\mu$
possesses an asymptotically stable single-round periodic orbit.

2$)$ The boundaries of domains $\triangle_k$ are bifurcation lines $L_k^+$ and $L_k^-$ corresponding to codimension
one bifurcations, saddle-node and period-doubling ones, respectively.

3$)$ As $k \to \infty$ the domains $\triangle_k$ accumulate to curve $L_h \cup \{L_+ \cap \{\mu_1 < 0\} \}$.

\end{thm}

Figure~\ref{fig:bifdiag} illustrates the statement of Theorem~\ref{thm:main}. Also a possible geometry of the first return map $T_k$
while varying parameters inside $\triangle_k$ is represented there. Vertical dashed lines
$\displaystyle \mu_2 = \mbox{const} \sim \frac{1}{k^2}$ are nominal boundaries such that the corresponding domains $\triangle_k$ 
do not cross them.
Two domains are highlighted at the figure: $RD$ (rescaling domain) and $D_{LS}$ (Lukjanov-Shilnikov domain) in which
the geometry of the first return map is considerably different. Namely, if $\mu \in \mbox{RD}$, the geometry is similar
to that observed during the emergence of the Smale horseshoe in various problems related to the study of homoclinic tangencies.
Here the rescaling method can be applied (see lemma~\ref{lem:resc}) which shows that the first return map can be
represented in the form of H\'enon map. When $\mu \in D_{LS}$ the geometry is different: the image of a stripe under the
first return map will have a form of a ``half-horseshoe'' (see figure~\ref{fig:bifdiag}). Such geometry was observed, in
particular, in \cite{LukSh78}. Generally speaking, domains $RD$ and $D_{LS}$ do not intersect but bifurcation curves can be
smoothly continued from $RD$ to $D_{LS}$. Note that domain $D_{LS}$ appears at the bifurcation diagram due to the well-known
``effect of a neighbourhood'' when orbits leave the neighbourhood under consideration $\Pi^+$ without bifurcating and therefore
the information about them gets lost. However, if the neighbourhood size increases or parameters decrease, this
effect disappears. Anyway, we can consider any sufficiently small neighbourhood of the origin $\mu = 0$ as $RD$.

\begin{rem}
In the domain $RD$ for any fixed $\mu_2 < 0$ in the corresponding family $f_{\mu_1}$, an infinite cascade of periodic sinks
is observed $($a sequence $\delta_k = \triangle_k \cap \{\mu_2 = \mbox{const}\})$ according to the famous theorem by
Gavrilov and Shilnikov {\rm \cite{GaS73}}. Analogous cascade of bifurcations will be observed also in family
$f_{\mu_2}$ for fixed $\mu_1 > -\varepsilon$. However, for $\mu \in D_{LS}$ in family $f_{\mu_2}$ a cascade
$\tilde\delta_k = \triangle_k \cap \{\mu_1 = \mbox{const}\}$ will be observed in accordance to the corresponding theorem
by Lukjanov and Shilnikov {\em \cite{LukSh78}}. Also note that only a finite number of periodic sinks will be observed in any
one-parametric family $f_{\mu_1}$ for fixed $\mu_2 > 0$.
\end{rem}


%

\section{Properties of maps $T_0$ and $T_1$}\label{sect:local}

It is well-known \cite{Luk79}, \cite{LukSh78} that $C^{r - 1}$-smooth coordinates exists in $U_0$ such that the local map $T_0$ can be written, for all sufficiently small $\mu$, in the following (finitely smooth) normal form
%
\begin{equation}
\label{eq:localnorm}
\left\{ \begin{array}{l}
\bar x = A(\mu) x + O(\|x\|^2 |y|) \\
\bar y = \mu_2 + y + y^2 + O(|y|^3),
\end{array}
\right.
\end{equation}
where $A(\mu)$ is an $(n\times n)$-matrix such that $\|A\| = \lambda < 1$.
Note that the second equation of (\ref{eq:localnorm})
does not depend on $x$-coordinates. In the coordinates
$(x, y)$ at $\mu = 0$ the fixed point $O$ lies in the origin and
its unstable $W^u$ and strong stable $W^{ss}$ manifolds are
locally straightened, their equations are $W^u:  \{x = 0\}$, $W^{ss}: \{y = 0\}$. For all small $\mu$, 
the invariant foliation $F^{ss}$ 
consists of $n$-dimensional planes $\{y = \mbox{const}\}$ and the axis $x = 0$ is also an invariant line.
For $\mu_2 < 0$  the local map $T_0$ possesses two fixed points, saddle $O_1$ and node $O_2$:
\begin{equation}
\label{eq:points}
O_1: \; \{ x = 0, y = \sqrt{-\mu_2} + O(|\mu_2|^{3/2})\}, \;
O_2: \; \{ x = 0, y = -\sqrt{-\mu_2} + O(|\mu_2|^{3/2})\}
\end{equation}
According to (\ref{eq:localnorm}) and (\ref{eq:points}) the local stable and unstable manifolds 
$W^s_{loc}(O_1)$ and $W^u_{loc}(O_1)$ have equations
$\{ y = \sqrt{-\mu_2} + O(|\mu_2|^{3/2})\}$ and 
$\{x = 0\}$, respectively.

By (\ref{eq:localnorm}), we can assume that $M^+ = (x^+,0)$ and $M^-= (0,y^-)$, 
where $x^+ > 0$ and $y^- > 0$. We note that points from $\Pi^+$ can reach $\Pi^-$ under iterations of $T_0$. 
In order to construct the corresponding map from $\Pi^+$ to $\Pi^-$ by the orbits of $T_0$, we discuss, 
first, geometrical properties of $T_0$ for various $\mu_2$.


%

\begin{figure}[ht]
\begin{center}
\epsfig{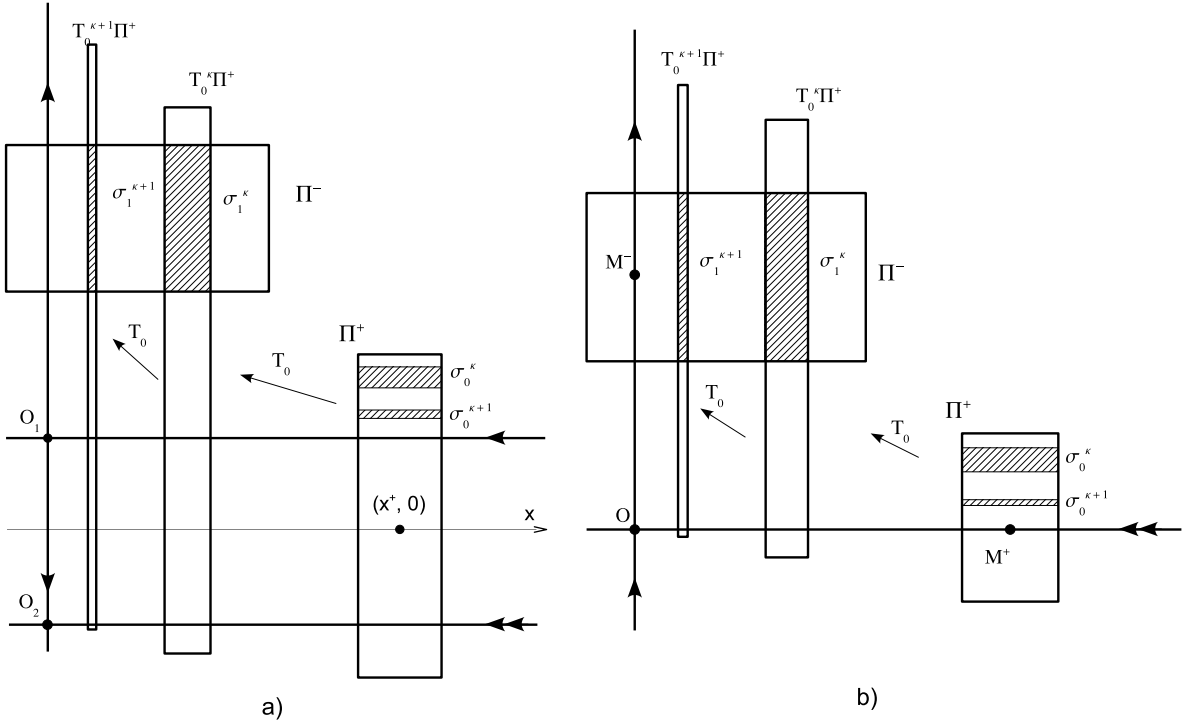}
\caption{{\footnotesize Geometry of the local map for a) $\mu_2 < 0$; b) $\mu_2 = 0$.}}
\label{fig:localfixed}
\end{center}
\end{figure}

If $\mu_2 \le 0$, then, evidently, there is such sufficiently large $\bar k$ that the images
$T_0^k \Pi^+$  for all $k\geq\bar k$
have non-empty intersections with  $\Pi^-$, see Fig.~\ref{fig:localfixed}. Correspondingly, in $\Pi^-$ there exist
a countable set of non-intersecting ``stripes'' $\sigma^1_k = T_0^k \Pi^+ \cap \Pi^-$ which accumulate
to $W^u_{loc}$  as $k \to \infty$.
The pre-images of these ``stripes'' are ``stripes''
$\sigma^0_k \equiv T_0^{-k}\sigma^1_k \subset \Pi^+$ which accumulate either to 
$W^{ss}_{loc}(O)$ at $\mu_2 = 0$ or to $W^s_{loc}(O_1)$ at $\mu_2 < 0$.

\begin{figure}[ht]
\begin{center}
\epsfig{file=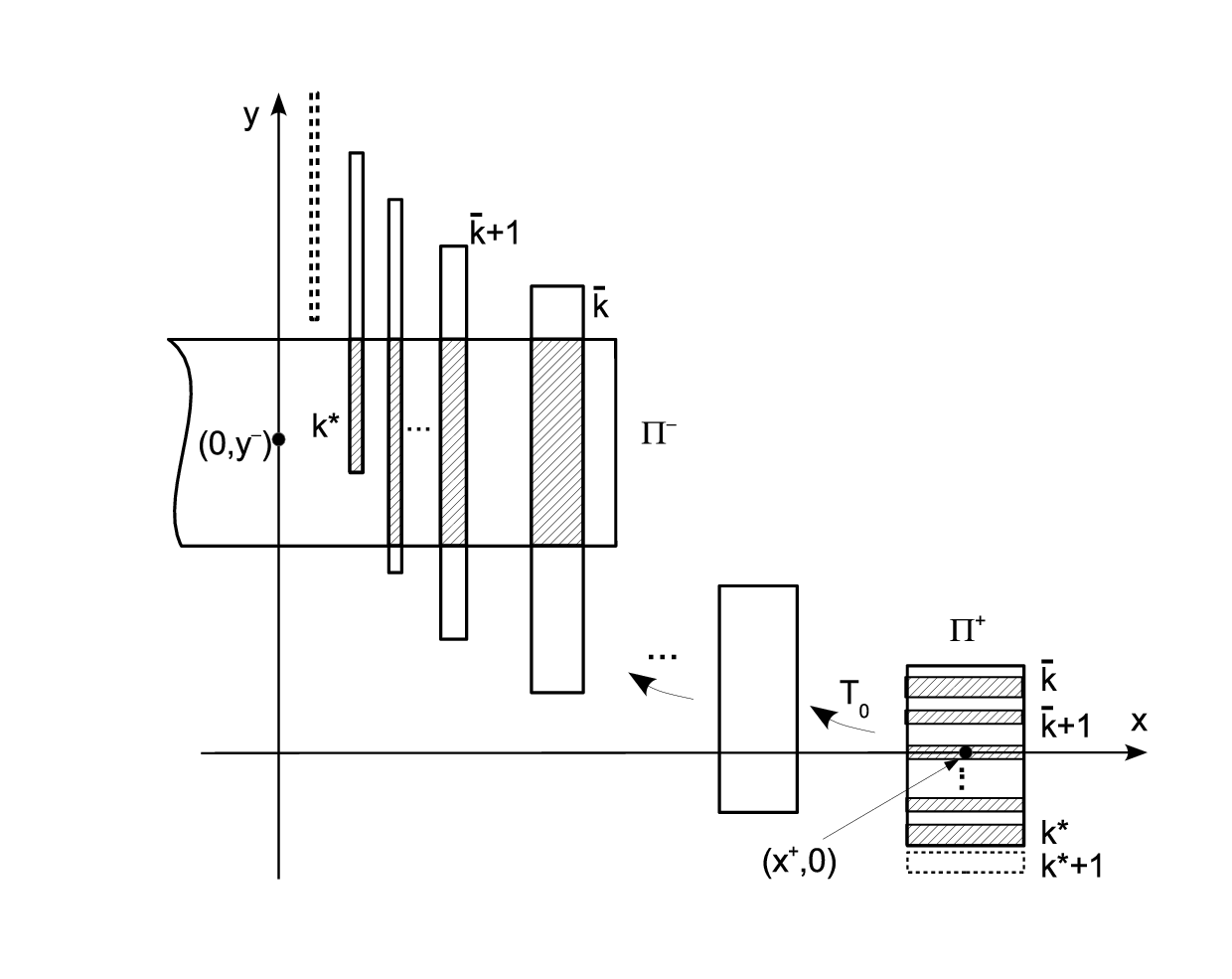,width=12cm}
\caption{{\footnotesize Geometry of the local map for $\mu_2 > 0$.}}
\label{fig:localpos}
\end{center}
\end{figure}

The situation is very different when $\mu_2 > 0$, see Figure~\ref{fig:localpos}, as any orbit leave $U_0$ for a finite number
of iterations of $T_0$. Consequently, only a finite number of stripes $\sigma_k^0$ such that
$T_0^k \sigma_k^0 \cap \Pi^- \neq \emptyset$ remains in $\Pi^+$. Thus, for any $\mu_2>0$ there exists
some natural $k^* = k^*(\mu_2) > \bar k$, $k^* \to \infty$ as $\mu_2 \to +0$, such that $T_0^k \Pi^+ \cap \Pi^- = \emptyset$ for all $k > k^*$.

Let $(x_i, y_i)$, $i = \overline{0, k}$ be a set of points in $U_0$ such that $(x_{i+1},y_{i+1}) = T_0(x_i,y_i)$.
One of the important peculiarities of local coordinates (\ref{eq:localnorm}) is that the relation
$(x_k, y_k) = T_0^k(x_0, y_0)$ allows a quite convenient representation in the so-called Shilnikov cross-form \cite{Sh67}. 
Thus, the following lemma holds.


%



\begin{lem}
\label{lem:local}
For any sufficiently small $\mu$, the  map $T_0^k: (x_0, y_0) \to (x_k, y_k)$ can be represented in the following
form:

\begin{equation}
\label{eq:loc}
\begin{array}{l}
x_k = A^k x_0  + o(\lambda^k) \xi_k(x_0, y_k, \mu)  \\

\displaystyle y_0 = \frac{\nu_k(\mu) y^- - \mu_2}{\nu_k(\mu) + y^-} +
 \frac{\theta_k(\mu) (1 + \rho_k(\mu))}{(\nu_k(\mu) + y^-)^2}(y_k - y^-) +
\theta_k(\mu) O((y_k - y^-)^2),

\end{array}
\end{equation}
%
%
where

\begin{equation}
\label{eq:loc14}
\nu_k(\mu) = \left\{
\begin{array}{ll}
\displaystyle \frac{\sqrt{-\mu_2}}{\tanh[k \sqrt{-\mu_2}]}, &  \mu_2 < 0\\

\displaystyle \frac{1}{k}, &  \mu_2 = 0 \\

\displaystyle \frac{\sqrt{\mu_2}}{\tan[k \sqrt{\mu_2}]}, &  \mu_2 > 0,
\end{array}
\right.
\end{equation}
\begin{equation}
\label{eq:loc21}
\theta_k(\mu) = \nu_k^2 + \mu_2,
\end{equation}
$|\rho_k(\mu)| << 1$, functions $\xi_k(x_0, y_k, \mu)$ are
uniformly bounded in $k$ along with all their derivatives up to the order $(r - 3)$, $\|x_k\|_{C^{r-2}} = O(\lambda^k)$
and $\|x_k\|_{C^{r-1}} \to 0 $ as $k \to \infty$.

\end{lem}

Note that an analogous result was obtained in \cite{Gor07a}, \cite{Gor07b} for the case $\mu_2\geq 0$. 
From the second equation of (\ref{eq:loc}), one can obtain the following
estimate for $k^*$ (see also \cite{Gor07a}, \cite{Gor07b}):

%
%
%
%
%
%
%
\begin{equation}
\label{eq:loc4}
\displaystyle k^{*} \sim \frac{\pi}{\sqrt{\mu_2}} - \frac{1}{\varepsilon}.
\end{equation}
%

Formulas (\ref{eq:loc14}) and (\ref{eq:loc21}) imply that for all sufficiently small $\mu$
$$
\theta_k = \left\{
\begin{array}{lc}
\displaystyle \frac{-\mu_2}{\sinh^2 k \sqrt{-\mu_2}},  & \mu_2 < 0\\
\displaystyle \frac{1}{k^2},  & \mu_2 = 0  \\
\displaystyle \frac{\mu_2}{\sin^2 k \sqrt{\mu_2}},  & \mu_2 > 0.\\
\end{array}\right.
$$

%
%
%

The global map $T_1: \Pi^- \to \Pi^+$ for all sufficiently small $\mu$ can be written
in the following form:
\begin{equation}
\label{eq:globgen}
\begin{array}{l}
\bar x_0 - x^+ = F(x_1, y_1 - y^-, \mu)\\
\bar y_0 = G(x_1, y_1 - y^-, \mu),
\end{array}
\end{equation}
where $(x_0, y_0) \in \Pi^+$, $(x_1, y_1) \in \Pi^-$;
functions $F$ and $G$ are $C^{r - 1}$-smooth and $F(0, 0, 0) = G(0, 0, 0) = 0$.
Condition {\bf B} means that $W^u(O)$ and $W^{ss}(O)$ at $\mu = 0$ have a quadratic tangency at the point $M^+(x^+, 0)$ that implies 
$G'_y(0, 0, 0) = 0$, $G''_{yy}(0, 0, 0) \neq 0$.

We can rewrite (\ref{eq:globgen}) as follows:
\begin{equation}
\label{eq:global}
\begin{array}{l}
\bar x_0 - x^+ = a x_1 + b (y_1 - y^-) + O(\|x_1\|^2 + \|x_1\|(y_1 - y^-) + (y_1 - y^-)^2 )\\
\bar y_0 = y^+(\mu) + c^{\top} x_1 + d(y_1 - y^-)^2 + O(\|x_1\|^2 + \|x_1\|(y_1 - y^-) + (y_1 - y^-)^3 ).
\end{array}
\end{equation}
where $y^+(0) = 0$ and the coefficients $y^+, y^-, x^+, d$ as well as $n$-dimensional vectors $b, c$ and $n \times n$ matrix $a$ depend smoothly on $\mu$.
The condition {\bf B} means that $d \neq 0$. 
Note that $d > 0$ (resp. $d < 0$) if the homoclinic tangency is 
from the saddle zone  (the node zone).


Without loss of generality, we set $\mu_1 \equiv y^+(\mu)$. Then, as it is seen from (\ref{eq:global}), $\mu_1$ is the distance
between $T_1W^u_{loc}(O)$ and $W^{ss}_{loc}(O)$ at $\mu_2 = 0$.

\subsection{Proof of Theorem~\ref{thm:bifdiag}}
By virtue of (\ref{eq:points}) the equation of $W^s_{loc}(O_1) \cap \Pi^+$ for $\mu_2 < 0$ is
$\{\|x - x^+\| < \varepsilon, \; y = \sqrt{-\mu_2} + O(|\mu_2|^{3/2})\}$. By (\ref{eq:global}), 
the equation of the curve $T_1(W^u_{loc}(O_1))$ is
$\{x_0 = x^+ + b(y_1 - y^-) + O((y_1 - y^-)^2), \; y_0 = \mu_1 + d (y_1 - y^-)^2 + O(|y_1 - y^-|^3), \;
|y_1 - y^-| \le \varepsilon\}$. Thus,
$T_1(W^u_{loc}(O_1))$ has a quadratic tangency with $W^u_{loc}(O_1)$ for $\mu_1 = \sqrt{-\mu_2} + O(|\mu_2|^{3/2})$ ---
this is the equation of the bifurcation curve $L_h$. By definition, $\mu_2 = 0$ is the
equation of the bifurcation curve $L^+$.

\section{The first return maps and the description of their bifurcations}


Using formulas (\ref{eq:loc}) and (\ref{eq:global}) we can easily construct (for all sufficiently large $k$ and small $\mu$) the
first return maps $T_k = T_1 T_0^k: \sigma^0_k \to \Pi^+$ in the local coordinates (\ref{eq:localnorm}).
However, the obtained formula, see (\ref{eq:resc2}), will be not very suitable for calculations 
because of its complexity. So that we will apply the so-called {\em rescaling method} \cite{TLY86,GST07} to this formula.
The main idea of this method is to perform smooth changes of coordinates and parameters (renormalization)
in order to bring the map to some standard form where small terms do not affect the dynamics.
The following lemma formalizes this approach to our case.

\begin{lem}
\label{lem:resc}
{\em (rescaling lemma)} Let $f_\mu$ be the family of diffeomorphisms defined above.
On the parameter plane $(\mu_1, \mu_2)$ there exist domain $RD$ $($Rescaling Domain$)$ that includes the origin, where,
using smooth changes of coordinates $(x_0, y_k) \to (X, Y)$ and parameters, for sufficiently large $k$,
the first return map $T_k$ can be brought to the following form:
\begin{equation}
\label{eq:parabola}
\begin{array}{l}
\bar X = b Y + o(1)_{k \to \infty}, \\
\bar Y = M - Y^2 + o(1)_{k \to \infty},
\end{array}
\end{equation}
where 
new coordinates $X=(X_1,...,X_n)$, $Y$ and parameter $M$
are defined in a ball $\|X, Y, M\| < S_k$, where $S_k \to +\infty$ as $k \to \infty$.
The $o(1)$-terms here denote some functions of $(X, Y, M)$ which tend to zero as $k \to \infty$
together with all their derivatives $($up to the order $(r - 2))$, and
\begin{equation}
\label{param}
M = \displaystyle -\frac{(y^-)^4}{d \theta_k^2}(\mu_1 - \nu_k + \ldots).
\end{equation}
\end{lem}

The proof of the rescaling lemma is given in section~\ref{sect:lemresc}.

Lemma~\ref{lem:resc} shows that the study of bifurcations of the first return map $T_k$ for large $k$
can be reduced to the study of the standard parabola map:
%
\begin{equation}
\label{eq:parabola1}
\bar Y = M - Y^2
\end{equation}
bifurcations of which are well-known. Namely, when $M \in (-1/4, 3/4)$ the map (\ref{eq:parabola1}) possesses a stable
fixed point which is born under the saddle-node bifurcation at $M = -1/4$ and loses stability via the period-doubling bifurcation
at $M = 3/4$. 
From formula (\ref{param}) we obtain that the equations of the bifurcation curves $L_k^+$
corresponding to a saddle-node bifurcation and $L_k^-$ corresponding to a period-doubling bifurcation 
have the following form:
\begin{equation}
\label{eq:bifdiag}
\begin{array}{l}
\displaystyle L_k^+: \;\; \mu_1 =
\nu_k(\mu_2) +
\frac{d}{4(y^-)^4} \theta_k^2(\mu_2) + \ldots \\
\\
\displaystyle L_k^-: \;\; \mu_1 =
\nu_k(\mu_2) -
\frac{3d}{4(y^-)^4} \theta_k^2(\mu_2) + \ldots \\
\end{array}
\end{equation}

Figure~\ref{fig:bifdiag} illustrates the bifurcation diagram of single-round periodic orbits 
in the case of tangency from the saddle zone ($d > 0$). In the case of tangency from the node zone ($d < 0$)
it will look quite similar with the only difference that curves $L_k^+$ and $L_k^-$ will change  places
($L_k^+$ will be the lower and $L_k^-$ the higher boundaries of  $\triangle_k$, see Figure~\ref{fig:tongue}).

Theorem~\ref{thm:main} is proved. It remains only to prove the main technical results --- lemmas~\ref{lem:local}
and~\ref{lem:resc}.

\section{Proof of lemma~\ref{lem:local}}
\label{sec:lemma1}

Consider an orbit $(x_i, y_i)$, $i = \overline{0, k}$, such that
$(x_0, y_0) \in \sigma_0^k \subset \Pi^+$, $(x_k, y_k) \in \sigma_1^k \subset \Pi^-$ and
$(x_{i + 1}, y_{i + 1}) = T_0 (x_i, y_i)$.


Provided that the second equation in (\ref{eq:localnorm}) does not depend on coordinates $x$, we will
consider it first. We will prove formula (\ref{eq:loc}) assuming that inequalities $|y_0| < \varepsilon$ and
$|y_k - y^-| < \varepsilon$ hold in $\Pi^+$ and $\Pi^-$. Here, without loss of generality, we assume that
$\Pi^+$ and $\Pi^-$ are sufficiently small rectangular neighbourhoods having the same fixed diameter $\varepsilon$.
We rewrite the second equation from (\ref{eq:localnorm}) as
\begin{equation}
\label{eq:loc1}
\bar y = y + g(y, \mu), \; \mbox{where} \; g(y, \mu) = \mu_2 + y^2 + O(y^3).
\end{equation}
It easy to
see that $g(y, \mu) > 0$ for $\mu_2 > 0$. If $\mu_2 \le 0$ we exclude from our consideration the domain 
of attraction of the saddle-node (for $\mu_2 = 0$) or the node (for $\mu_2 < 0$) as it
obviously does not contain points, iterations of which (with respect to $T_0$) reach $\Pi^-$.
In other words we consider only a sub-domain of $U^-$ where $y > \sqrt{-\mu_2} + O(|\mu_2|^{3/2})$ i.e.
above $W^{ss}(O)$ or $W^s(O_1)$ respectively. For these values of $y$ for all small $\mu$ function $g(y, \mu)$
is strictly positive so that sequence $\{ y_0, \ldots, y_k \}$ is monotonically increasing.

In order to find a formula for $y_0$ for all small $\mu_2$ we consider
the following integral from $y$ to $\bar y$:
\begin{equation}
\label{eq:loc17}
\displaystyle \int\limits_{y}^{y + g(y, \mu)} \frac{dy}{\mu_2 + y^2} = 1 + \phi(y, \mu),
\end{equation}
where
\begin{equation}
\label{eq:loc19}
\phi(y, \mu) =
\left\{\begin{array}{ll}
\displaystyle \frac{1}{\sqrt{\mu_2}} \arctan \frac {\sqrt{\mu_2}}{1 + y + O(y^3)} - 1, &
\mbox{if} \; \mu_2 > 0,\\

\displaystyle  \frac {O(y)}{1 + y + O(y^2)}, & \mbox{if} \; \mu_2 = 0,\\

\displaystyle \frac{1}{2\sqrt{-\mu_2}} \log \left[\frac{1 + \sqrt{-\mu_2} + y + O(y^3)}
{1 - \sqrt{-\mu_2} + y + O(y^3)}\right] - 1, & \mbox{if} \; \mu_2 < 0.\\
\end{array}\right.
\end{equation}

Formula (\ref{eq:loc17}) can be written for each $y_0, y_1, \ldots, y_{k - 1}$ by taking $y = y_{i}$ and $y + g(y, \mu) = y_{i + 1}$.
Adding together all these integrals for $i = 0, 1, \ldots, k - 1$
we obtain the following implicit dependency of $y_0$ on $y_k$:
\begin{equation}
\label{eq:loc18}
\displaystyle \int\limits_{y_0}^{y_k} \frac{dy}{\mu_2 + y^2} = k + \eta(y_k, \mu), \; \mbox{where} \;
\eta(y_k, \mu) = \sum\limits_{i = 0}^{k - 1} \phi(y_i(y_k, \mu), \mu),
\end{equation}

Note that functions $\phi(y, \mu)$ and $\eta(y, \mu)$ are $C^{r - 1}$-smooth. The required smoothness of function
$\phi(y, \mu)$ follows from the fact that the integral in (\ref{eq:loc19}) is taken from a smooth function.
Then, we have $\eta(y, \mu) \in C^{r - 1}$ as it is a finite sum of smooth functions $\phi(y_i(y_k, \mu), \mu)$.

Introduce the following value:
\begin{equation}
\label{eq:loc6}
\displaystyle \zeta_k(\mu) = \sum \limits_{i = 0}^{k - 1} \frac{1} {g(y_i(y^-, \mu), \mu)}.
\end{equation}

\begin{lem}
\label{lem:zero}
Function $\eta(y_k, \mu)$ is uniformly bounded together with its derivatives with respect to $y_k$ up to order $(r - 1)$.
For the derivatives of $\eta(y_k, \mu)$ with respect to the parameters the following uniform in $k$ estimates hold:
\begin{equation}
\label{eq:loc8}
\displaystyle \left| \frac{\partial^{k_1 + k_2 + k_3}\eta(y_k, \mu)}
{\partial \mu_1^{k_1} \partial \mu_2^{k_2} \partial y_k^{k_3}}
\right|  \le \mbox{const} \cdot \zeta_k^{k_1 + k_2}(\mu)
\end{equation}
\end{lem}

{\bf Proof.}

Derivatives of $\phi(y, \mu)$ with respect to $y_k$ are obtained via an explicit differentiation of
integral (\ref{eq:loc18}) with respect to $y_k$. They are easily estimated (see \cite{ST00}) as follows:
%
%
%
%
%
%
\begin{equation}
\label{eq:loc9}
\displaystyle \left| \frac{\partial^{m}\phi(y_j(y_k, \mu), \mu)}{\partial y_k^{m}}
\right|  \le \mbox{const} \cdot g(y_j, \mu).
\end{equation}
Thus for the derivatives of $\eta(y_k, \mu)$ we obtain:
\begin{equation}
\label{eq:loc10}
\begin{array}{c}
\displaystyle \left| \frac{\partial^{m}\eta(y_k, \mu)}{\partial y_k^{m}} \right|  \le
\sum\limits_{i = 0}^{k - 1}\left| \frac{\partial^{m}\phi(y_j(y_k, \mu), \mu)}{\partial y_k^{m}} \right|
 \le \mbox{const} \cdot \sum\limits_{i = 0}^{k - 1} g(y_j, \mu) = \\

= \mbox{const}\cdot \sum\limits_{i = 0}^{k - 1} (y_{i + 1} - y_i) =  \mbox{const} \cdot (y_k - y_0)
 \le \mbox{const} \cdot (y^- + 2\varepsilon).
\end{array}
\end{equation}

Also, in analogy to \cite{ST00} it can be shown that
%
\begin{equation}
\label{eq:loc11}
\displaystyle \left| \frac{\partial^{k_1 + k_2 + k_3}y_j(y_k, \mu)}
{\partial \mu_1^{k_1} \partial \mu_2^{k_2} \partial y_k^{k_3}}
\right|  \le \mbox{const} \cdot g(y_j, \mu) \zeta_k^{k_1 + k_2}(\mu).
\end{equation}
Estimates of the same kind are valid for the derivatives of $\phi(y_i(y_k, \mu), \mu)$.
Adding together these inequalities we finally obtain estimate (\ref{eq:loc8}). The lemma is proved.

Lemma~\ref{lem:zero} allows us to represent formula (\ref{eq:loc18}) in the following form:
\begin{equation}
\label{eq:loc12}
\displaystyle \int\limits_{y_0}^{y_k} \frac{dy}{\mu_2 + y^2} = k + \eta(y_k, \mu) = k (1 + \delta_k) + O(y_k - y^-),
\end{equation}
where $\delta_k \to 0$ as $k \to \infty$. Following this we rewrite the local map
(\ref{eq:loc18})as follows:
\begin{equation}
\label{eq:loc13}
 y_0 = \left\{
\begin{array}{ll}
\displaystyle
  \frac{y_k \sqrt{\mu_2} - \mu_2 \cdot \tan\left[ k \sqrt{\mu_2} (1 + \delta_k)
 + \sqrt{\mu_2}O(y_k - y^-)\right]}
{y_k \cdot \tan\left[ k \sqrt{\mu_2} (1 + \delta_k) + \sqrt{\mu_2}O(y_k - y^-)\right] + \sqrt{\mu_2}},
&   \mu_2 > 0\\
\\
\displaystyle \frac{y_k} {y_k  ( k (1 + \delta_k) + O(y_k - y^-)) + 1}, &   \mu_2 = 0\\
\\
\displaystyle \frac{y_k \sqrt{-\mu_2} - \mu_2 \cdot \tanh\left[ k \sqrt{-\mu_2} (1 +
 \delta_k) + \sqrt{-\mu_2}O(y_k - y^-)\right]}
{y_k \cdot \tanh\left[ k \sqrt{-\mu_2} (1 + \delta_k) +
\sqrt{-\mu_2}O(y_k - y^-)\right] + \sqrt{-\mu_2}}, &   \mu_2 < 0\\
\end{array}\right.
\end{equation}

Expanding formula (\ref{eq:loc13}) into the Taylor series near
point $y_k = y^-$ we obtain formula (\ref{eq:loc}) with $\nu_k$ given by formula (\ref{eq:loc14}) and
$|\rho_k| \sim g^2(y^-, \mu) << 1$.

Now formula for $x_k$ from (\ref{eq:loc}) can be obtained as
a solution of the corresponding boundary value problem \cite{GS93}, \cite{book} for which coordinates
$y_i$, $i = 0, 1, \ldots, k$ are already known. Lemma~\ref{lem:local} is proved.

\section{Proof of Lemma~\ref{lem:resc}}
\label{sect:lemresc}

By virtue of formulas (\ref{eq:loc}) and (\ref{eq:global}) we can write first return map $T_k$ in the following form:
\begin{equation}
\label{eq:resc1}
\begin{array}{l}
\bar x - x^+ = a A^k x + o(\lambda^k)\xi_k(x, y, \mu) + b (y - y^-) +
 O(\lambda^{2k} \|x\|^2  +\lambda^k \|x\||y - y^-| + \\ \qquad + (y - y^-)^2 )\\

\displaystyle
\frac{\nu_k y^- - \mu_2}{\nu_k + y^-} +
 \frac{\theta_k (1 + \rho_k)}{(\nu_k + y^-)^2} (\bar y - y^-) +
\theta_k O((\bar y - y^-)^2) = \mu_1 + c A^k x + \\
 \qquad + o(\lambda^k) \xi_k(x, y, \mu) + d(y - y^-)^2 + O(\lambda^{2k}\|x\|^2 + \lambda^k \|x\| |y - y^-|
 + \\ \qquad + |y - y^-|^3 ).
\end{array}
\end{equation}
We perform a shift of coordinates $x_{new} = x - x^+ + \psi_{1k}, \; y_{new} = y - y^- + \psi_{2k}$,
where $\psi_{1k}, \psi_{2k} = O(\lambda^k)$,
to eliminate the constant term in the first equation and the linear term in $y_{new}$ in the
second equation of (\ref{eq:resc1}). The map is rewritten as:
\begin{equation}
\label{eq:resc2}
\begin{array}{l}
\bar x = a A^k x + b y +
o(\lambda^{k})O(\|x\|^2) + \lambda^k O(\|x\| |y|) + O(y^2)\\

\displaystyle
 \frac{\theta_k(1 + \rho_k)}{(\nu_k + y^-)^2} \bar y +
 \theta_k O(\bar y^2) = \\
\qquad  = M_1 + c A^k x + dy^2 + o(\lambda^{k})O(\|x\|^2) + \lambda^k O(\|x\| |y|) + O(|y|^3),
\end{array}
\end{equation}
where $\displaystyle M_1 = \mu_1 - \nu_k(\mu) + \ldots$ and the dots stand for higher order terms.
Now rescale the coordinates using the following formulas:
\begin{equation}
\label{eq:resc3}
x = \alpha  X, \; y = \alpha Y, \\
\end{equation}
where 
%
\begin{equation}
\label{eq:resc4}
\alpha = \displaystyle -\frac{\theta_k(1 + \rho_k)}{d (y^-)^2}.
\end{equation}
After this, at those $k$ for which $\theta_k$ are asymptotically small,
map (\ref{eq:resc2}) can be represented in the form (\ref{eq:parabola})
and formula (\ref{param}) is valid for $M$. We denote as $RD$ the domain of the parameter plane where
$\theta_k \to 0$ as $k \to \infty$. It is easy to check that this domain includes the origin but is bounded from below
by some curve of the form $\mu_1 = -\varepsilon + O(\sqrt{\mu_2})$ which passes through point $(\mu_1 = -\varepsilon, \mu_2 = 0)$.
We also note that this curve is a conditional boundary of domain $D_{LS}$, where the rescaling-method does not work as
$\theta_k$ is not infinitesimal here for $k \to \infty$.
Lemma~\ref{lem:resc} is proved.

\subsection*{Acknowledgements.}
The authors are grateful to Prof. D.V. Turaev for very useful discussions.
The paper was supported by grant 14-41-00044 of the RSF,
grants of RFBR, projects No.13-01-00589,
13-01-97028--povolzhje and 14-01-00344, EPSRC Mathematics Platform grant EP/I019111/1
and the Leverhulme Trust grant RPG-279.


\begin{thebibliography}{40}


\bibitem{Sh67}
L.P. Shilnikov. On a Poincare-Birkhoff problem. {\em Math. USSR Sb.}, 1967, vol. 3, No.3, pp. 353--371.
%
\bibitem{GTS91}
S.V. Gonchenko, D.V. Turaev, L.P. Shilnikov,
On models with non-rough Poincare homoclinic curves (in Russian), {\em Dokl. Akad. Nauk}, 1991, vol. 320,
pp. 269--272.
%
\bibitem{GST93} S.V. Gonchenko, L.P. Shil'nikov, D.V. Turaev,
On models with non-rough Poincare homoclinic curves, {\em Physica D.}, 1993,
vol. 62, pp. 1--14.
%
\bibitem{GTS99}
S.V. Gonchenko, D.V. Turaev, L.P. Shilnikov,
Homoclinic tangencies of arbitrary order in Newhouse domains (in Russian),
{\em Proc. of Int. Conf. dedicated to L.S. Pontryagin 90th anniversary}, 1999, vol. 67, pp. 69--128.
%
\bibitem{GST07} S.V. Gonchenko, L.P. Shilnikov, D.V. Turaev,
Homoclinic tangencies of arbitrarily high orders in
conservative and dissipative two-dimensional maps, {\em Nonlinearity}, 2007, vol. 20,
pp. 241--275.
%
\bibitem{GaS73} N.K. Gavrilov, L.P. Shilnikov, On three-dimensional dynamical
systems close to systems with a structurally unstable homoclinic
curve, Part I. {\em Math. USSR Sbornik}, 1972, vol. 17,
467--485; Part II. {\em ibid}, 1973, vol. 19, pp. 139--156.
%
\bibitem{LukSh78}
Lukjanov V.I., Shilnikov L.P., On some bifurcations of dynamical systems with homoclinic structures,
{\em DAN SSSR}, 1978, vol. 243, no. 1, pp. 26--29.
%
\bibitem{G83}
{\em S.V.Gonchenko} On stable periodic motions in systems close to
a system with a nontransversal homoclinic curve. Russian Math.
Notes. 1983. V.33. No.5. P.384--389.
%
\bibitem{AfSh83} Afraimovich V.S., Shilnikov L.P.,
Invariant tori, their breakdown and stochasticity,
{\em Amer. Math. Soc. Transl.}, 1991, vol. 149, pp. 201--211.
%
\bibitem{NPT} S. Newhouse, J. Palis, F. Takens, Bifurcations and stability of families of diffeomorphisms.
{\em IHES Publ. Math.}, 1983, vol. 57, 5--71.
%
\bibitem{ST00} Shilnikov, L.P., Turaev, D.V., A new simple bifurcation of a periodic orbit of ``blue sky
Catastrophe'' type, Methods of qualitative theory of differential equations and related topics ,
Amer. Math. Soc. Transl. Ser. 2, 200, Amer. Math. Soc., Providence, RI, 2000, pp. 165-188.
%
\bibitem{Gor07a}
Gordeeva O.V., Lukjanov V.I.,
On bifurcations of dynamical systems of codimension two having a non-rough homoclinic
structure of a saddle-node,
{\em Vestnik Nizhegorodskogo universiteta}, 2007, no. 2, pp. 175--180.
%
\bibitem{Gor07b}
Gorgeeva O.V., Lukjanov V.I.,
Certain bifurcations of limit sets in a neighbourhood of a non-rough homoclinic
structure with a non-degenerated periodic motion,
{\em Nelineinyi mir}, 2007, vol.5, no. 1--2, pp. 95--100.
%
\bibitem{Luk82}
V.I. Lukjanov, Bifurcations of dynamical systems with a saddle-node separatrix loop,
{\em Differentsialnyje uravnenija}, 1982, vol. 18, no. 9, pp. 1493--1506; English transl.,
{\em Differential equations}, 1983, vol. 18, pp. 1049--1059.
%
\bibitem{TLY86} Tedeschini-Lalli L., Yorke J.A.
How often do simple dynamical processes have infinitely many coexisting sinks? {\em Commun.Math.Phys.} 1986. V.106. P.635-657.
%
\bibitem{Luk79}
Lukjanov V.I., On existence of smooth invariant foliations in a neighbourhood of
some non-rough fixed points of a diffeomorphism,
{\em Differential and integral equations,
Mezhvuz. sb.}, Gorky, 1979, vol. 3.
%
\bibitem{HPS}
M.W. Hirsch, C.C. Pugh, M. Shub, {\em Lecture Notes in Math.}, v.583,
Springer-Verlag, Berlin, 1977.
%
\bibitem{book} L.P. Shilnikov, A.L. Shilnikov, D.V. Turaev, L.O. Chua, {\em Methods
of qualitative theory in nonlinear dynamics. Part~I}, World Scientific, 1998.
%
%
\bibitem{GS93}
S.V. Gonchenko, L.P. Shilnikov, On moduli of systems with a structurally unstable homoclinic
Poincare curve, {\em Russian Acad. Sci. Izv. Math.}, 1993, vol. 41, no. 3, pp. 417--445.








\end{thebibliography}
\end{document}